\documentclass[12pt,reqno]{amsart}

\usepackage{amsmath}
\usepackage{amssymb}
\usepackage{amsthm}
\usepackage[english]{babel}

\newtheorem{theorem}{Theorem}

\newtheorem{corollary}{Corollary}
\newcommand{\beqa}{\begin{eqnarray}}
\newcommand{\beqan}{\begin{eqnarray*}}
\newcommand{\eeqa}{\end{eqnarray}}
\newcommand{\eeqan}{\end{eqnarray*}}
\def\beq#1\eeq{\begin{equation}#1\end{equation}}

\def\P{\mathbf P }

\def\g{\gamma }

\def\be{\beta}

\def\la{\lambda}

\def\me{^{-1}}

 \def\na{\,\, {\raise.4pt\hbox{$\shortmid$}}{\hskip-2.0pt\to}\, \, }

\def\={\overset{ \text{\rm def} }=}

\def\ffrac{\frac}
\def\4{\kern1pt}

\def\bgr#1{\4\bigr#1}
\def\bgl#1{\bigl#1\4}
\def\BR{\bf R}

\begin{document}

\title[Rare events]
{Rare events and Poisson point
processes}

\author[F.~G\"otze]{Friedrich G\"otze}
\author[A.Yu. Zaitsev]{Andrei Yu. Zaitsev}

\email{goetze@math.uni-bielefeld.de}
\address{Fakult\"at f\"ur Mathematik,\newline\indent
Universit\"at Bielefeld, Postfach 100131,\newline\indent D-33501 Bielefeld,
Germany\bigskip}
\email{zaitsev@pdmi.ras.ru}
\address{St.~Petersburg Department of Steklov Mathematical Institute
\newline\indent
Fontanka 27, St.~Petersburg 191023, Russia\newline\indent
and St.Petersburg State University, 7/9 Universitetskaya nab., St. Petersburg,
199034 Russia}

\begin{abstract}
The aim of the present work is to show that the results obtained earlier 
on the approximation of distributions of sums of independent terms by 
the accompanying compound Poisson laws may be interpreted as rather sharp 
quantitative estimates for the closeness between the sample containing 
independent observations of rare events and the Poisson point process 
which is obtained after a Poissonization of the initial sample.
\end{abstract}

\keywords {sums of independent random variables, rare events, Poisson point processes, 
concentration functions, inequalities}

\subjclass {Primary 60F05; secondary 60E15, 60G50}

\thanks{The authors were supported by the SPbGU-DFG grant 6.65.37.2017.
The second author was supported by grant RFBR 16-01-00367 and by the Program of the Presidium
of the Russian Academy of Sciences No. 01 'Fundamental Mathematics and its Applications' under grant
PRAS-18-01.}

\maketitle

The aim of the present work is to show that the results obtained earlier 
on the approximation of distributions of sums of independent terms by the accompanying 
compound Poisson infinitely divisible laws may be interpreted as rather sharp quantitative 
estimates for the closeness between the sample containing independent observations 
of rare events and the Poisson point process which is obtained after a Poissonization of the initial sample.

Let us first introduce some notation. Let $\frak F_d$ denote the set
of probability distributions defined on the Borel $\sigma$-field of
subsets of the Euclidean space~${\BR}^d$ and \,${\mathcal L}(\xi)\in\frak F_d$ \,is the
distribution of a $d$-dimensional random vector~$\xi$. For $F\in\frak F\=\frak F_1$,  the
concentration function is defined by \;${Q(F,\,b)=\sup_x
F\bgl\{[x,\,x+b]\bgr\}}$, \, $b\ge0$. For $F\in\frak F_d$, we
denote the corresponding  distribution functions
by~$F(x)=F\{(-\infty,x_1]\times\cdots\times(-\infty,x_d]\}$, $x=(x_1,\ldots,x_d)\in{\BR}^d$, 
and the uniform Kolmogorov distance by
$$ \rho (F,H)=\sup_x \bgl|F(x)-H(x)\bigr|. $$ By~\,$c$ \,we
denote absolute positive constants. Note that constants \,$c$ \,can be different  
in different (or even in the same) formulas. For
some positive quantities $a$ and $b$ writing $a\asymp b$ means
that $a\le c\,b$ and $b\le c\,a$.

Let $X_1, X_2,\dots, X_n$ be independent not
identically distributed elements of a measurable space $(\frak
X,\mathcal S)$ and $f:\frak X\to\bold R^d$ be a  Borel mapping.
 Let $F_i={\mathcal L}(f(X_i))$, $i=1,2,\dots,n$,
be the distributions of $f(X_i)$. Then the sum $$ S=f(X_1)+
f(X_2)+\dots+ f(X_n)$$ has the distribution $\prod_{i=1}^nF_i$
(products and powers of measures will be understood in the
convolution sense: \;${FH=F*H}$, \,$H^m=H^{m*}$, \,$H^0=E\=E_0$,
where $E_a$ \,is the distribution concentrated at a point $a\in
{\BR}^d $). A natural approximating infinitely divisible distribution
for $\prod_{i=1}^nF_i$ is the accompanying compound Poisson
distribution $\prod_{i=1}^n e(F_i)$, where $$
e(H)\=e^{-1}\sum_{m=0}^\infty
 \ffrac{H^m}{m!},\quad H\in\frak F_d,
$$ and, more generally, $$ e(\alpha\4 H)\=e^{-\alpha}\sum_{m=0}^\infty
 \ffrac{\alpha^m\4H^m}{m!}, \quad \alpha>0.
$$ It is easy to see that $\prod_{i=1}^n e(F_i)$ is the
distribution of \begin{equation}T=\sum_{i=1}^n\sum_{j=1}^{\nu_i}f(X_{i,j}), \label{(135}\end{equation}
where $X_{i,j}$ and $\nu_i$, $i=1,\dots,n$, $j=1, 2, \dots$, are
independent in aggregate random elements of the space $\frak X$
and random variables respectively with ${\mathcal L}(X_{i,j})={\mathcal L}(X_{i})$ and
${\mathcal L}(\nu_i)=e(E_1)$. Clearly, $ e(E_1)$ is the Poisson distribution
with mean 1. Thus, the sum $T$ is defined similarly to $S$, but
the initial sample $\bold X=(X_1, X_2,\dots, X_n)$ is replaced by
its Poissonized version $\bold Y=\bgl\{X_{i,j}:$ $i=1,\dots,n$,
$j=1, 2, \dots,\nu_i\bgr\}$. Poissonization of the sample is known
as one of the most powerful tools in studying empirical processes.
The random set $\bold Y$ may be considered as a realization of the Poisson point
process on the space $\frak X$ with intensity measure
$\sum_{i=1}^n {\mathcal L}(X_{i})$. The important property of the Poisson point process 
is the space independence: for any pairwise disjoint sets $A_1,\ldots,A_m \subset\frak X$, 
the random sets $\mathbf Y\cap A_1,\ldots,\mathbf Y\cap A_m \subset\frak X$ are independent in aggregate. 
As a consequence, the Poisson point
process $\bold Y$ allows for better approximation than the process $\bold X$. 
One can use the independence property since  the theory of independent objects is much more elaborated.

In this paper, we consider the problem of approximation of the sample 
by the Poisson point process which is obtained after a Poissonization 
of the initial sample in the case where the sample is obtained by observation of rare events. 
The situation is the following.
The set $\frak X$ is represented as the union of two disjoint measurable sets: $\frak X= \frak
X_1\cup\frak X_2$, with \,$\frak X_1,\,\frak X_2\in \mathcal S$, $\frak
X_1\cap\frak X_2=\varnothing$. We say that the $j$-th  rare event occurs if $X_j\in\frak X_2$. 
Respectively, it does not occur if $X_j\in\frak X_1$.

To formulate the
results we need some additional notation.

Let $f:\frak X\to\bold R^d$ be a Borel mapping defined above 
and $F_i={\mathcal L}(f(X_i))$, $i=1,2,\dots,n$. Then
 distributions $F_i\in\frak F_d$ can be represented as mixtures \begin{equation}
F_i=(1-p_i)\,U_i+p_i\,V_i, \label{(11}\end{equation} where $U_i, V_i\in\frak F_d$ 
are conditional distributions of $f(X_i)$ given
$X_i\in\frak X_1$ and $X_i\in\frak X_2$ respectively, that is
\begin{equation} 0\le p_i=\mathbf{P}\big\{X_i\in\frak X_2\big\}
=1-\mathbf{P}\big\{X_i\in\frak X_1\big\}\le1. \label{(12}\end{equation}
Below the  $V_{i}$ are arbitrary distributions. We deal with rare events while the quantity
\begin{equation}  p=\max_{1\le i\le
n}p_i\label{(13}\end{equation}
 is small. In other words, this is the case  if our rare events are sufficiently rare.
Let
\begin{equation}
a_i=\int_{\bold R^d} x\,U_i\{dx\},\quad i=1,2,\ldots, n,\label{(14}\end{equation}
and, for $d=1$, $$ {\left|\bf a\right|_2^2}=\sum_{i=1}^n
a_i^2,\quad \left|\bold a\right|_\infty=\max\limits_{1\le i \le n}|a_i|, $$ $$
\sigma_i^2=(1-p_i)\int_{-\infty}^{\infty}(x-a_i)^2\,U_i\{dx\},\quad
B^2=\sum_{i=1}^n\sigma_i^2, $$ Denote \begin{equation}
H_1=\prod_{i=1}^nF_i,\quad H_2=\prod_{i=1}^n e(F_i),\quad
H_3=\prod_{i=1}^n E_{a_i}e(F_iE_{-a_i}). \label{(15}\end{equation}

Usually, a good approximation of the distribution~$H_1$ is given by the distribution~$H_3$. 
But while estimating the closeness of  the sample $\bold X$ to the Poisson point process $\bold Y$, 
we are interested in the closeness of distributions~$H_1$ and~$H_2$. 
Thus, we need to estimate the distances between $H_1$ and~$H_2$
using
that $H_2$ is the distribution of $$
T=\sum_{i=1}^n\sum_{j=1}^{\nu_i}f(X_{i,j}), $$ where $f(X_{i,j})$ and
$\nu_i$, $i=1,\dots,n$, $j=1, 2, \dots$, are random vectors and variables with ${\mathcal L}(f(X_{i,j}))=F_i$,
${\mathcal L}(\nu_i)=e(E_1)$, which are independent in
aggregate. Moreover, $ H_3={\mathcal L}(T-\Delta)$, with
$\Delta=\sum_{i=1}^n a_i\,(\nu_i-1)$.

There exists a lot of results on the compound Poisson approximation of
the distributions $\prod_{i=1}^nF_i$, see Arak and Zaitsev \cite{2},
 Barbour and Chryssaphinou~\cite{3}, 
Zaitsev \cite{14}, Roos~\cite{Roos}, \v Cekanavi\v cius~\cite{Ch16} and
the bibliography therein. However, most of these results
require that distributions $F_i$ are appropriately centered. In the
general case a better approximation can be obtained for the
accompanying compound Poisson approximation of
$H_1$ by $H_3$ using
centering constants $a_i$.

In this paper, we consider estimates of classical distances between 
the distributions $H_1$ and~$H_2$ with remainder terms having the additive 
summand of the form~$ c\,p$ for $d=1$ or~$ c(d)\,p$ for $d\ge1$. The most natural 
result of such a type is given in the following theorem.

\begin{theorem}\label{Theorem0} Let the conditions above be satisfied for $d\ge1$ and $f(x)=0$, 
for all $x\in\frak X_1$. Then
\begin{equation} \rho\big(H_1,\,H_2 \big) \le c(d)\,p.  \label{76}\end{equation}
\end{theorem}
Theorem~\ref{Theorem0} is a direct consequence of a result of Zaitsev~\cite{z89} 
which provides inequality~\eqref{76} in the case, where \begin{equation}
U_i=E, \qquad
i=1,2,\dots,n, \label{789}\end{equation} and the  $V_{i}$ are arbitrary distributions. 
For $d=1$, this statement was obtained earlier in~\cite{z83}.  
This was an improvement of a result of Le Cam~\cite{9} who has proved one-dimensional assertion 
with $p^{1/3}$ instead of~$p$.

Approximation of the sample by a Poisson point process under the conditions of Theorem~\ref{Theorem0} 
was discussed in Zaitsev~\cite{14}, see also Hipp~\cite{Hipp} and Roos~\cite{Roos}.

Inequality~\eqref{76} shows that different samples with rare events of uniformly small~$p_i$ 
are close one to another if they have the same (or close) measures $\sum_{i=1}^n p_i\,V_i$. 
Indeed, in this case they have the same (or close) approximating Poisson point processes $\bold Y$. 
In the case of identical distributions~$V_i$, inequality~\eqref{76} can be essentially sharpened and improved. 
Even the distance in variation $d_{TV} \big(H_1,\,H_2 \big)$ can be estimated by~$ c\,p$ 
(Prokhorov~\cite{p53} for the case of identical~$p_i$ and Le Cam~\cite{L60} 
for the case of arbitrary~$p_i$) and, moreover, by $\sum_{i=1}^{n}p_{i}^{2}\Big/\max\{1, \sum_{i=1}^{n}p_{i}\}$ 
(Barbour and Hall~\cite{4}). Thus, for identical~$V_i$, the statement on the closeness of our processes 
can be essentially strengthened. The condition that the rare events are identically distributed is, however, 
seldomly satisfied. Rare events are similar to extreme incidents; each of them is unique and has his own 
individual distribution. Thus, inequality~\eqref{76} can be useful, for example, in insurance theory 
to estimate the probabilities that the cumulative influence of risk factors $f(X_j)$ will 
not exceed fixed critical values $x_j$ (see, e.g., Hipp~\cite{Hipp}).

In Theorems~\ref{Theorem1}--\ref{45} below, the remainder terms of the form~$ c\,p$ or~$ c(d)\,p$ 
are related to Theorem~\ref{Theorem0} which is used in the proofs.

The importance of centering in the accompanying approximation was
noted by Le Cam~\cite{9}. He has shown that accompanying laws with
centering (see distribution $H_3$) provide the rate of
approximation of the form $c\,n^{-1/3}$ for distributions of sums of $n$ i.i.d.\
random variables with an absolute constant $c$ without any
assumptions on the distribution of summands. The same rate was
obtained earlier by Kolmogorov \cite{8} for the approximation by
infinitely divisible distributions. Le Cam's result shows that
accompanying laws are sufficient to prove such a rate. Later Arak \cite{1}
obtained the optimal rate $c\,n^{-2/3}$ for infinitely divisible
approximation (see Arak and Zaitsev \cite{2} for the history of the
problem). Arak's infinitely divisible approximation is much more
complicated than the approximation by accompanying laws. Ibragimov
and Presman [6] have shown that the rate $c\,n^{-1/3}$ is optimal
for the accompanying approximations (see Arak and Zaitsev (\cite{2},
inequality (1.3), p.~181)).

Le Cam \cite{9} considered not only i.i.d.\  summands, he obtained bounds
for the accompanying compound Poisson approximation for
distributions of sums of independent non-identically distributed
random variables---with and without centering.

In Example 1 of G\"otze and Zaitsev \cite{GZ04}, the case where $F_1=\dots=F_n=F$,
$\int|x|^3\,F\{dx\}<\infty$, $a=\int x\,F\{dx\}\ne0$ and
$\sigma^2=\int (x-a)^2\,F\{dx\}>0$ was considered.  Using the Berry--Esseen bound in the CLT,
it was shown that
$\rho\big(H_1,\,H_2 \big)\ge
c\,\min\{1,\,a^2/\sigma^2\}$, for sufficiently large $n\ge n_0$.
In Example 2 of G\"otze and Zaitsev \cite{GZ04}, it was shown that in  the degenerate case, for example, if
$F=E_1$, we have $H_1=E_n$ and  $H_2=e(nE_1)$, the Poisson law with parameter $n$. Clearly,
$\rho(E_n,e(nE_1))\ge c$.

Thus, the approximation without centering is
not always successful. To ensure the validity of such an
approximation, one needs additional assumptions such as zero mean,
symmetry or a large atom at zero, see Zaitsev \cite{14} and the
bibliography therein.

Note, however, that these examples are rarely appearing in our scheme of rare events.
It is much more plausible that the values of the function~$f(x)$ are much larger for $x\in \frak X_2$ 
than for $x\in \frak X_1$. Moreover, the tails of distributions $V_i$ usually will be heavier 
than those of~$U_i$. It is well-known that the concentration functions of convolutions of distributions 
with heavy tails decrease faster than those of distributions with light tails. 
Thus, we may expect that the concentration functions of distributions $H_1,H_2,H_3$ 
admit good bounds for the remainder terms in Theorems~\ref{Theorem2}--\ref{Theorem5} below.

Le Cam~\cite{9} considered the rate of approximation by accompanying
laws without centering proving the following result (see Le Cam
\cite{10}, Theorem 2, p. 431).

\begin{theorem}\label{Theorem1} Let the above conditions be satisfied for $d=1$ and,  for
some $\tau>0$, \begin{equation}
U_i\bgl\{[-\tau,\,\tau]\bgr\}=1, \qquad
i=1,2,\dots,n. \label{(1)}\end{equation}
 Then
$$ \rho\big(H_1,\,H_2 \big) \le c\,\bgl (p^{1/3}+
\big(\big(1+\tau^{-2}{\left|\bf a\right|_2^2}\big) D^{-2}(\tau)\big)^{1/3}\bigr),
\eqno(2) $$ where $$
D^{2}(\tau)=\sum_{i=1}^n\int_{-\infty}^{\infty}\min\bgl\{1,
x^2\tau^{-2}\bgr\} F_i\{dx\}. $$
\end{theorem}

Applying Theorem \ref{Theorem1} in the case $a_i=0$, $i=1,2,\dots,n$, we obtain
another result of Le Cam \cite{9} which contains the rate of
accompanying approximation with centering, see Le Cam (\cite{10},
Theorem 1, p.~429):

\begin{corollary}\label{Corollary1} Let the conditions of Theorem\/ $\ref{Theorem1}$  be
satisfied with $\bf a=0$.
 Then
$$ \rho\big(H_1,\,H_2 \big) \le c\,\bgl
(p^{1/3}+ D^{-2/3}(\tau)\bigr). \eqno(3) $$
\end{corollary}

For the proof it suffices to notice that the distributions
$F_iE_{-a_i}$ satisfy the conditions of Theorem 1 with $a_i=0$ and
with replacing $\tau$ by $2\4\tau$.

The statement of \ref{Corollary1}  was slightly improved by Ibragimov
and Presman \cite{6}. Zaitsev \cite{12, 13} proved Theorem \ref{Theorem2} below which
sharpened and generalized the statement of Corollary~\ref{Corollary1}. To
formulate this theorem we need some additional notation.

Following Katz \cite{7} and Petrov \cite{11}, we introduce the class $\mathcal
G$ of real valued functions $g(\4\cdot\4)$ on $\bold R$ satisfying the conditions

a) $g(\4\cdot\4)$ is a non-negative even function, which is
strictly positive for $x\ne 0$ and does not decrease for $x\ge0$.

b) the function $x/g(x)$ is non-decreasing for $x>0$.

\begin{theorem}[\rm Zaitsev \cite{13}]\label{Theorem2}
 Let the conditions of Theorem\/ $\ref{Theorem1}$  be satisfied without condition~$\eqref{(1)}$.
 Let $g\in\mathcal G$. Denote
$$ \be_i=\be_i(g)=(1-p_i)\int_{-\infty}^{\infty}(x-a_i)^2\,
g(x-a_i)\,U_i\{dx\}<\infty, $$ $$ \be=\be(g)=\sum_{i=1}^n\be_i(g).
$$ Define
$\lambda=\lambda(g)=\min\big\{B,\ffrac{\be(g)} {B\,g(B)}\big\}$, if
$B^2>0$,  and $\lambda=0$ otherwise. Let $B^2>0$. Then \begin{equation} \rho\big(H_1,\,H_3 \big)\le c\,\bgl (p+\min\big\{Q(H_1,\la),Q(H_3,\la)\bigr\}\bigr)\asymp
 p+\ffrac{\la} {B}\,Q\Big(\prod_{i=1}^ne\big( p_i\,V_i\,E_{-a_i}\bigr),B\Big).
\label{(4)}\end{equation}
If $B^2=0$, then  $\rho\big(H_1,\,H_3
\big)\le c\4p$.
\end{theorem}

It can be shown that, under condition \eqref{(1)}, $Q(H_3,\la) \le c\,(p+D\me(\tau))$, see \cite{GZ04}. 
Therefore, Theorem~\ref{Theorem2} is sharper with respect to order
than Corollary~\ref{Corollary1}. The main result of Zaitsev \cite{12} is a particular
case of Theorem~\ref{Theorem2} with $g(x)\equiv |x|$. The optimality of 
the statement of Theorem~\ref{Theorem2} was analyzed in \cite[Examples~3 and~4]{GZ04}.

\medskip

Theorem \ref{Theorem2} is the main tool in the proof of the main results of
the paper~\cite{GZ04}, Theorems~\ref{Theorem3}--\ref{Theorem6}, which improve Theorem \ref{Theorem1}, see a discussion in~\cite{GZ04}.

\begin{theorem}\label{Theorem3} Let the conditions of Theorem\/ $\ref{Theorem2}$  be
satisfied. Then $$ \rho\big(H_1,\,H_2 \big) \le c\,\bgl (p+
\P\bgl\{\left|\Delta\right|>\g\bgr\} +\min_{1\le
k\le3}Q(H_k,\g)\bigr), $$ for any $\g\ge\lambda$, where
$\Delta=\sum_{i=1}^n a_i\,(\nu_i-1)$ and $\nu_i$ are i.i.d.\  Poisson
with mean $1$.
\end{theorem}

A particular case of Theorem \ref{Theorem3} occurring when $a_i=0$,
$i=1,2,\dots,n$, implies inequality~\eqref{(4)}. Theorems \ref{Theorem4}--\ref{Theorem6}
below are simple consecutive consequences of Theorem \ref{Theorem3}.

\begin{theorem}\label{Theorem4} Let the conditions of Theorem\/ $\ref{Theorem2}$  be
satisfied. Then, for any $\varkappa\ge\la$, $$ \rho\big(H_1,\,H_2
\big) \le c\,\bgl (p+(1+{\left|\bf a\right|_2} \varkappa\me\sqrt{\delta}+\left|\bold a\right|_\infty
\varkappa\me\delta)\,q \bigr),  $$ where
$q=\min\limits_{1\le k\le3}Q(H_k,\varkappa)$  and $\delta=\log
\bgl(1+\varkappa\, q\me{\left|\bf a\right|_2}\me\bigr)$.
\end{theorem}

\begin{theorem}\label{Theorem5} Let the conditions of Theorem\/ $\ref{Theorem1}$  be
satisfied. Then $$ \rho\big(H_1,\,H_2 \big) \le c\,\bgl (p+(1+{\left|\bf a\right|_2}
\tau\me\sqrt{s}+\left|\bold a\right|_\infty \tau\me s)\,Q\, \bigr),  $$ where
$Q=\min\limits_{1\le k\le3}Q(H_k,\tau)$ and $s=\log \bgl(1+\tau\4
Q\me{\left|\bf a\right|_2}\me\bigr)$.
\end{theorem}

\begin{theorem}\label{Theorem6} Let the conditions of Theorem\/ $\ref{Theorem1}$  be
satisfied. Then $$ \rho\big(H_1,\,H_2 \big) \le c\,\bgl (p+(1+{\left|\bf a\right|_2}
\tau\me\sqrt{r}+\left|\bold a\right|_\infty \tau\me r)\,D\me(\tau) \bigr), $$ where
$r=\log \bgl(1+\tau\,D(\tau)\,{\left|\bf a\right|_2}\me\bigr)$.
\end{theorem}

Theorem \ref{Theorem6} is better with respect to order than Theorem \ref{Theorem1} (clearly,
$\left|\bold a\right|_\infty\le{\left|\bf a\right|_2}$). Theorems \ref{Theorem3}--\ref{Theorem5} 
provide sharper bounds.

Finally, we mention a result on strong approximation of sums of independent random vectors 
by infinitely divisible distributions.  Theorem  \ref{hhthm4}
 below is a consequence of the main result of Zaitsev \cite{z89}, Lemma~A of Berkes and Philipp~\cite{hh3}, 
 and the Strassen--Dudley theorem  (see Dudley~\cite{hh9}). 
 The norms in the space  $\mathbf{R}^{d}$ are denoted 
 $\left\| \,\cdot\,\right\|_2$ and $\left\| \,\cdot\,\right\|_\infty$.
\medskip

\begin{theorem}\label{hhthm4}  Let conditions~$\eqref{(11}$--$\eqref{(15}$ be satisfied for $d\ge1$, and let, 
for some $\tau\ge0$,
\[
 U_{i}\left\{
\left\{ x\in \mathbf{R}^{d}:\left\| x\right\|_2 \le \tau \right\}
\right\} =1, \quad i=1,2,\ldots,n,
\]
and the  $V_{i}\in\frak F_d$ are arbitrary distributions. Then, for each fixed $\lambda >0$, 
the random vectors  \,$S$, $T$ and \,$\Delta $ can be constructed on the same probability space so that
\[
\mathcal{L}(S )=H_1,\qquad\mathcal{L}(T )=H_2,\qquad \mathcal{L}(T-\Delta
)=H_3,
\]
where
$\Delta=\sum_{i=1}^n(\nu_i-1)\,a_i$ and $\nu_i$ are i.i.d.\  Poisson
with mean $1$, and
\begin{equation}
\mathbf{P}\big\{ \Vert S -T+\Delta \Vert_2 >\lambda \big\} \le
c(d)\,\Big( p+\exp \Big( -\frac{\lambda
}{c(d)\,\tau }\Big) \Big) +\sum_{i=1}^{n}p_{i}^{2}.
\label{infdiv}
\end{equation}
Moreover, for all $x\in\mathbf{R}^{d} $ and $\lambda >0$,
\begin{equation}
H_1(x) \le H_3(x+\lambda\cdot\mathbf{1})+
c(d)\,\Big( p+\exp \Big( -\frac{\lambda
}{c(d)\,\tau }\Big) \Big),
\label{20}
\end{equation}
\begin{equation}
H_3(x) \le H_1(x+\lambda\cdot\mathbf{1})+
c(d)\,\Big( p+\exp \Big( -\frac{\lambda
}{c(d)\,\tau }\Big) \Big),
\label{21}
\end{equation}
where $\mathbf{1}=(1,1,\ldots,1)\in\mathbf{R}^{d} $ and $c(d)$ depends only on $d$. If the distributions
$V_{i}$ are identical, then the term \,$\sum_{i=1}^{n}p_{i}^{2}$ \,in
\eqref{infdiv} can be dropped.\end{theorem}
\bigskip

Theorem~\ref{hhthm4} implies the following assertions about the closeness of distributions $H_1$ and~$ H_2$.

\begin{theorem}\label{45}Let the conditions of Theorem\/ $\ref{hhthm4}$  be
satisfied. Then
\begin{equation}
\mathbf{P}\big\{ \Vert S -T \Vert_2 >2\,\lambda \big\} \le
c(d)\,\Big( p+\exp \Big( -\frac{\lambda
}{c(d)\,\tau }\Big) \Big) +\sum_{i=1}^{n}p_{i}^{2}+\mathbf{P}\big\{ \Vert \Delta \|_2 \ge\lambda \big\}.
\label{50}
\end{equation}
If the distributions
$V_{i}$ are identical, then the term \,$\sum_{i=1}^{n}p_{i}^{2}$ \,in
\eqref{50} can be removed. Moreover, for all $x\in\mathbf{R}^{d} $ and $\lambda >0$,
\begin{equation}
H_1(x) \le H_2(x+2\,\lambda\cdot\mathbf{1})+
c(d)\,\Big( p+\exp \Big( -\frac{\lambda
}{c(d)\,\tau }\Big) \Big)+\mathbf{P}\big\{ \Vert \Delta \|_\infty \ge\lambda \big\},
\label{89}
\end{equation}
\begin{equation}
H_2(x) \le H_1(x+2\,\lambda\cdot\mathbf{1})+
c(d)\,\Big( p+\exp \Big( -\frac{\lambda
}{c(d)\,\tau }\Big) \Big)+\mathbf{P}\big\{ \Vert \Delta \|_\infty \ge\lambda \big\}.
\label{67}
\end{equation}
\end{theorem}

Choosing an appropriate $\lambda$, we see that inequalities~\eqref{50} and~\eqref{89}--\eqref{67} 
imply bounds for the Prokhorov aand L\'evy distances respectively.

    Using Bernstein's
inequality, we proved in~\cite{GZ04} that, for $d=1$, \begin{equation}
\P\bgl\{\left|\Delta\right|\ge\g\bgr\} \le 2\,\max\bgl\{e^{-\g^2/4{\left|\bf a\right|_2^2}
},\;e^{-\g/4\left|\bold a\right|_\infty }\bgr\}, \label{677}
\end{equation} for any $\g\ge0$. This inequality can be used for the estimation of the right-hand sides of inequalities~\eqref{50}--\eqref{67}. For $d>1$, we can apply inequality~\eqref{677} coordinatewise.

\end{document}